\documentclass[10pt,a4paper,twoside,leqno]{article}
\usepackage{amsfonts,amssymb}
\usepackage{eufrak}
\usepackage{amscd}

\newtheorem{defn}{Definition}

\vspace{0.5cm}

 4
\font\ebf=cmbx8
\font\erm=cmr8

\usepackage{graphicx}

\begin{document}

\thispagestyle{empty}

\noindent {\bf Cobweb  posets as noncommutative prefabs }

\vspace{0.7cm} {\it A. Krzysztof Kwa\'sniewski}

\vspace{0.2cm}

{\erm High School of Mathematics and Applied Informatics}

{\erm  Kamienna 17, PL-15-021 Bia\l ystok, Poland}

\vspace{0.2cm}

\noindent {\ebf Summary}

{\small A class of new type graded infinite posets with minimal
element is introduced. These  so called cobweb posets proposed
recently by the present author constitute a wide range of new
noncommutative and nonassociative prefab combinatorial schemes`
examples with characteristic graded sub-posets as primes. These
schemes are defined here via relaxing commutativity and
associativity requirements imposed on the composition in prefabs
by the fathers of this fertile concept. The construction and the
very first basic properties of cobweb prefabs are disclosed. An
another new type prefab example with single valued commutative and
associative composition is provided. "En passant" though not by
accident - we discover  new combinatorial interpretation of all
classical $F-nomial$ coefficients hence specifically incidence
coefficients of reduced incidence algebras of full binomial type
are given a new cobweb combinatorial interpretation also.

\vspace{0.2cm}

 AMS Classification Numbers: 05C20, 11C08, 17B56 .

 \vspace{0.7cm}

\section{Introduction} The concept of prefab (with associative and commutative
composition) was introduced in [1], see also [2,3]. Here we shall
deliver a class   of similar combinatorial structure of new type
based on the so called cobweb posets.  For the sake of
completeness we recall in Section 1. the definition of a cobweb
poset as well as a combinatorial interpretation of its
characteristic binomial-type coefficients  (for example-
fibonomial ones) [4,5].\\
In Section 2. after relaxing associativity and commutativity
requirements imposed on the composition in prefabs by the authors
of this concept [1] we observe that the vast family of all cobweb
posets becomes by construction a new type of $nonassociative$
noncommutative prefabs` subclass. The very first basic properties
of these cobweb prefabs are shown up. As a result a class of new
type of graded infinite posets with minimal element are employed
here as an enveloping framework for the completely new class of
combinatorial prefab structures with noncommutative and
$nonassociative$ composition (synthesis) of its objects since now
on called $\textbf{prefabiants}$. Cobweb infinite posets $P$ are
designated uniquely by any cobweb admissible sequence of integers
$F=\{n_F\}_{n\geq 0}$ and are by construction endowed with
self-similarity property. Namely at each graded level vertex a
family of infinite cobweb sub-posets isomorphic to $P$ may be rooted.\\
The number of finite characteristic sub-posets (\textrm{prime
$prefabiants$}) of  $P$ in between levels of graded Hasse digraph
of $P$ is given by $F-nomial$ coefficients. These include:
incidence coefficients such as binomial or  $q$-Gaussian ones for
finite geometries or fibonomial coefficients [4,5] which are not
incidence coefficients. Here these $F-nomial$ numbers are
introduced also via $c_2)$ axiom in the Definition 1 from [1]:

\begin{equation}
|a\odot b|= \frac{f(a\odot b)}{f(a)f(b)}= \left(
\begin{array}{c} n\\k\end{array}\right)_{F}.
\end{equation}
We notice with emphasis and not only occasionally,  that $c_2)$
\textbf{axiom} in Definition 1 from [1] \textbf{is equivalent to
the fundamental Theorem} 1 from [1]. More then that - in Section 3
we shall see that all objects from Equation 1 gain specific
uniform combinatorial interpretation within the class of cobweb
prefab combinatorial scheme - by construction.

\vspace{2mm}

\section{Cobweb posets - presentation and their combinatorial interpretation}

\vspace{1mm} Given any sequence $\{F_n\}_{n\geq 0}$ of nonzero
reals one defines  its corresponding  binomial-like $F-nomial$
coefficients in the spirit of  Ward`s Calculus of sequences
[6](reals may be replaced for example by any field of
characteristic zero) as follows

\begin{defn}
$$
\left( \begin{array}{c} n\\k\end{array}
\right)_{F}=\frac{F_{n}!}{F_{k}!F_{n-k}!}\equiv
\frac{n_{F}^{\underline{k}}}{k_{F}!},\quad n_{F}\equiv F_{n}\neq
0, n\geq 0 $$ where we make an analogy driven identifications in
the spirit of  Ward`s Calculus of sequences  $(0_F\equiv0)$:
$$
n_{F}!\equiv n_{F}(n-1)_{F}(n-2)_{F}(n-3)_{F}\ldots 2_{F}1_{F};$$
$$0_{F}!=1;\quad n_{F}^{\underline{k}}=n_{F}(n-1)_{F}\ldots (n-k+1)_{F}. $$
\end{defn}
This is just the adaptation of the notation for the purpose
Fibonomial Calculus case (see Example 2.1 in [7]). \vspace{2mm}
Given any such sequence $\{F_n\}_{n\geq 0}\equiv\{n_F\}_{n\geq 0}$
of now  \textbf{nonzero integers } we define following [4,5] the
partially ordered graded infinite set  $P$ - called afterwards a
cobweb poset - as follows. Its vertices are labelled  by pairs of
coordinates: ${\langle i , j \rangle} \in {N \times N_0}$ where
$N_0$ denotes the nonnegative integers. Vertices show up in layers
("generations") of $N \times N_0$ grid along the recurrently
emerging subsequent $s-th$ levels $\Phi_s$ where $s\in N_0$  i.e.

\begin{defn}
$$\Phi_s =\{\langle j, s\rangle 1\leq j \leq s_F\}, {s\in N_0}. $$
\end{defn}
We shall refer to $\Phi_s$  as to  the set of vertices at the
$s-th$ level. The population of the  $k-th$ level ("generation" )
counts  $k_F$  different member vertices for $k>0$ and one for
$k=0$. \vspace{2mm} Here down a disposal of vertices on $\Phi_k$
levels is visualized for the case of Fibonacci sequence (the
subtlety of $F_0=0$ is manageable) \vspace{3mm}

$---\Uparrow-----\Uparrow----up --Fibonacci---stairs--\star--k-th-level$\\

$---- and ----- so ---- on ---- up    --- \Uparrow ----------$\\
$\star \star \star \star \star \star \star \star \star \star \star \star \star \star \star \star \star \star\star \star \star \star \star \star \star \star \star \star \star \star \star \star \star \star \star \star \star \star \star \star \star \star \star \star --\star \star \star \star\star10-th-level$\\
$\star \star \star \star \star \star \star \star \star \star \star \star \star \star \star \star \star \star \star \star \star \star \star \star \star \star \star \star \star \star \star\star\star\star----------- 9-th-level$\\
$\star \star \star \star \star \star \star \star \star \star \star \star \star \star \star \star \star \star \star \star \star------------------8-th-level$\\
$\star \star \star \star \star \star \star \star \star \star \star \star \star -------------------------7-th-level$\\
$\star \star \star \star \star \star \star \star-----------------------------6-th-level$\\
$\star \star \star \star \star ---------------------------------5-th-level$\\
$\star \star \star ---------------------------------- 4-th-level$\\
$\star \star -----------------------------------3-rd-level  $ \\
$\star ------------------------------------ 2-nd-level$\\
$\star ----------------------------------- 1-st-level$\\
$\star ----------------------------------- 0-th-level$\\
 \vspace{2mm}
   \textbf{Figure 1. The $s-th$ levels in $ N\times N_0 $}
   , $N_0$ - nonnegative integers

\vspace{2mm}

\noindent Accompanying the set $E$ of edges to the set $V$ of
vertices - we obtain the Hasse diagram where here down ${p,q,s}\in
N_0 $. (Convention: Edges stay for arrows directed - say -
upwards) Namely:

\begin{defn}
$$P=\langle V,E\rangle,\quad V=\bigcup_{0\leq p}\Phi_p ,\quad E
=\{\langle\langle j , p\rangle ,\langle q ,(p+1) \rangle
\rangle\}\bigcup\{\langle\langle 1 , 0\rangle ,\langle 1 ,1
\rangle \rangle\},$$ \quad where $1 \leq j \leq {p_F} , 1\leq q
\leq {(p+1)_F}.$
\end{defn}

\begin{defn}
The finite  cobweb sub-poset $P_m = \bigcup_{0\leq s\leq m}\Phi_s$
is called the prime cobweb poset.
\end{defn}
In reference [3,4] a partially ordered infinite set $P$ was
introduced  via  descriptive picture of  its Hasse diagram.
Indeed, we may picture out the partially ordered infinite set $P$
from the  Definition $3$ with help of the sub-poset $P_{m}$ ({\it
rooted at $F_{0}$ level of the poset}) to be continued then ad
infinitum in now obvious way as seen from the figures  $Fig.1-
Fig. 5$ of $P_{m}$ cobweb posets below. These look like the
Fibonacci rabbits` way generated tree with a specific
``cobweb''[4,5,8]. This is an example of acyclic directed graphs
(DAG) [9] cobweb subclass.

\begin{center}
\includegraphics[width=75mm]{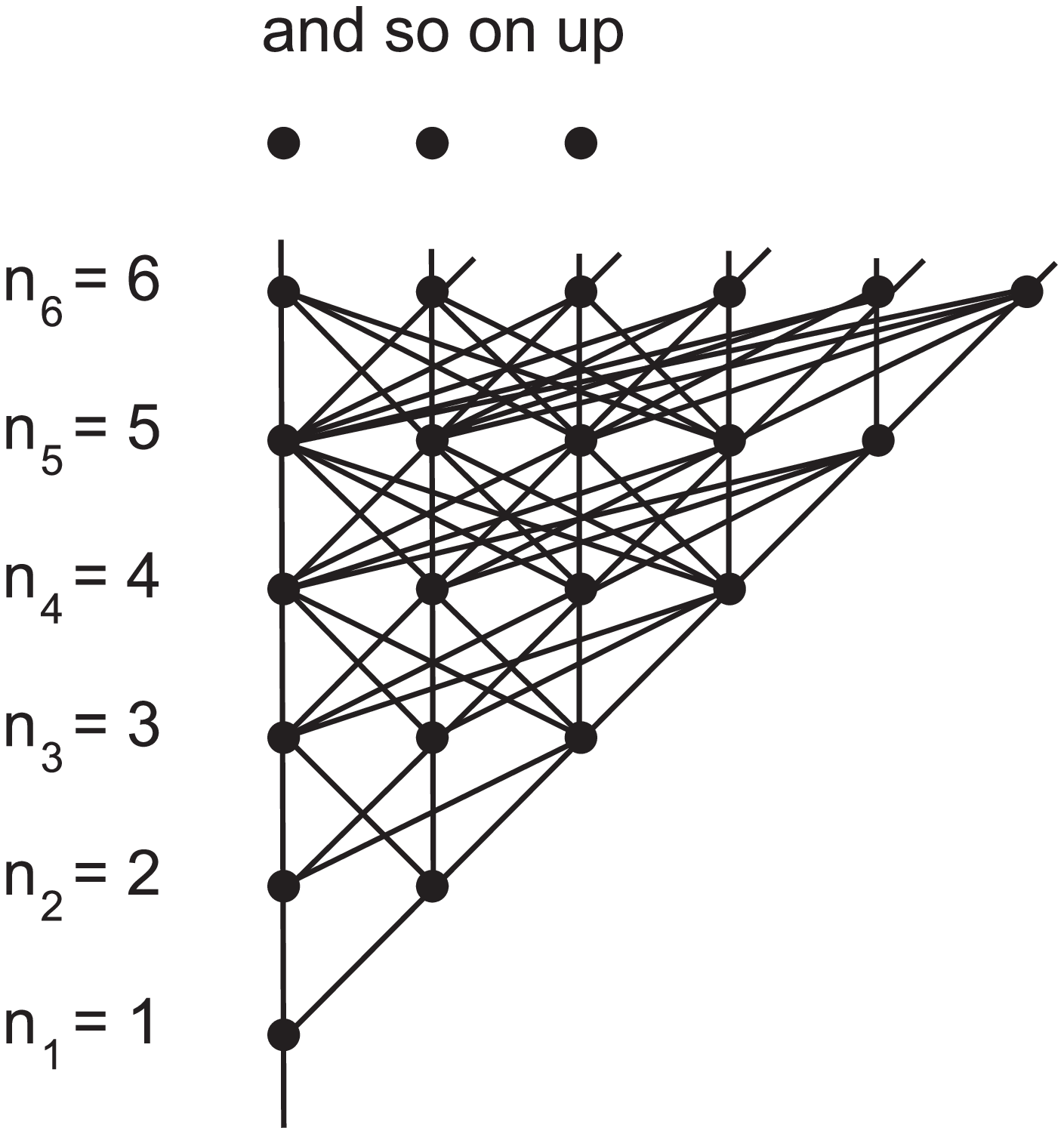}

\vspace{2mm}

\noindent {\small Fig.1. Display of Natural numbers` cobweb
poset.}
\end{center}

\vspace{2mm}

\begin{center}

\includegraphics[width=75mm]{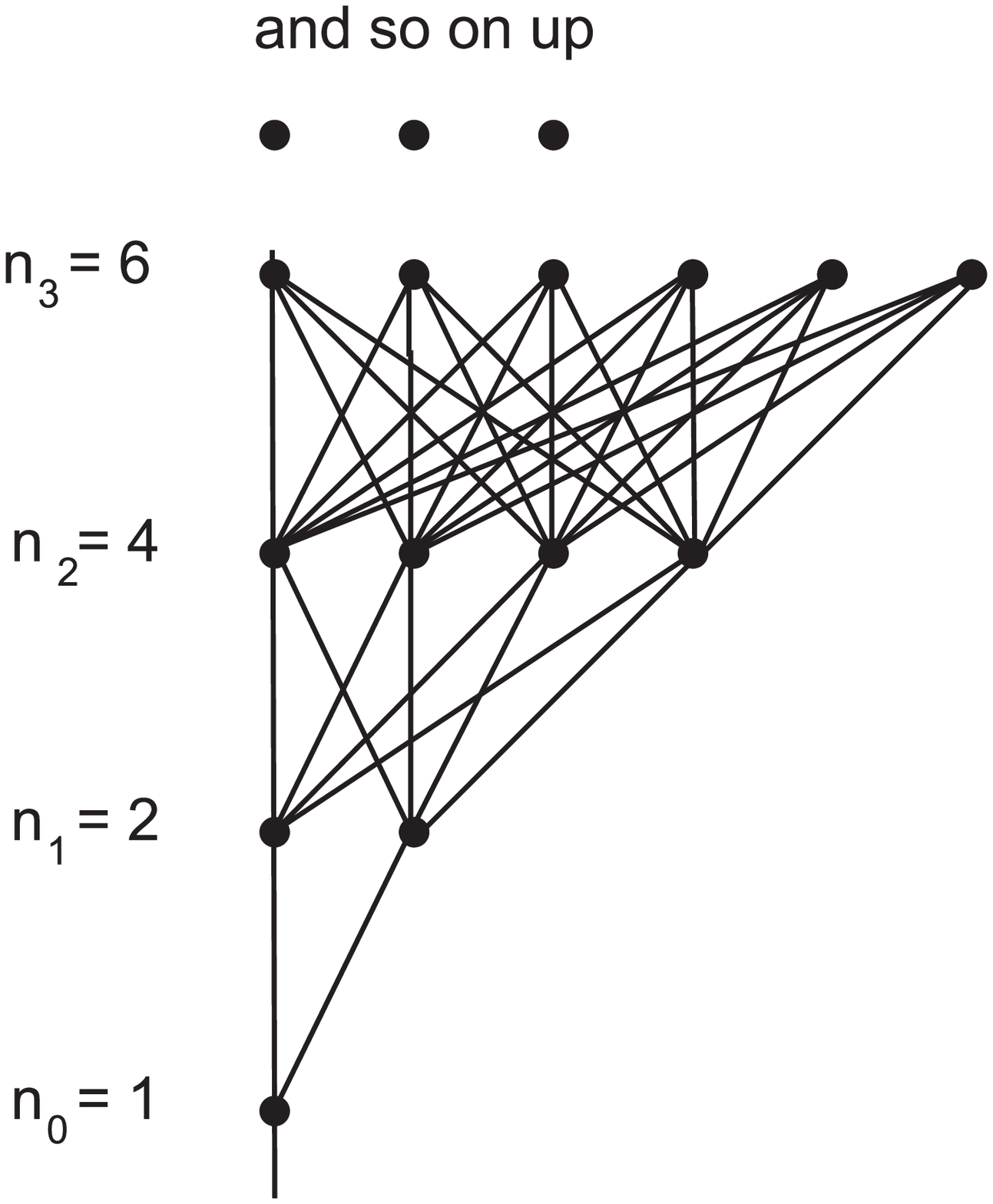}

\noindent {\small Fig.2. Display of Even Natural numbers` cobweb
poset.}

\end{center}

\vspace{2mm}

\vspace{2mm}

\begin{center}

\includegraphics[width=75mm]{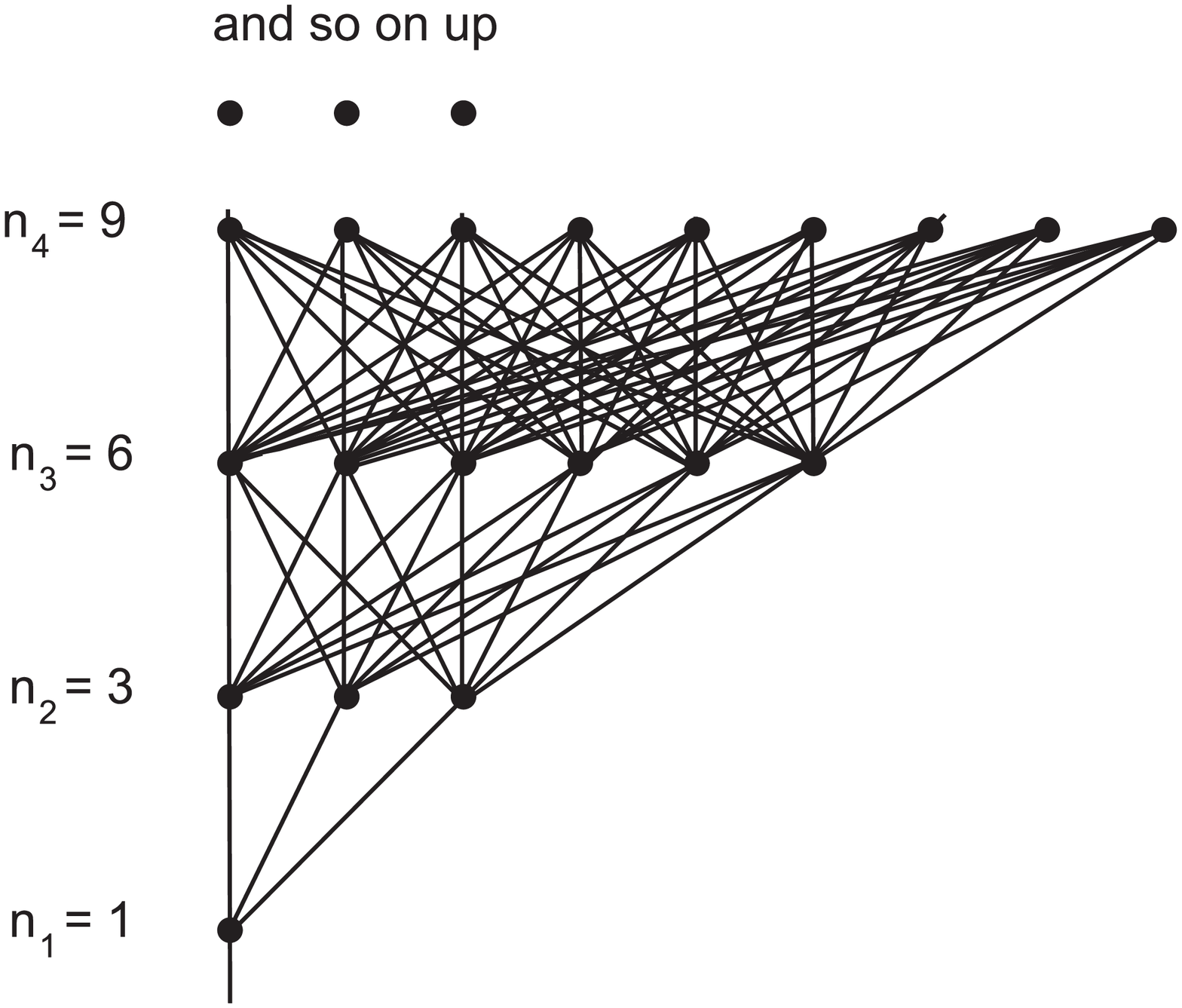}

\vspace{2mm}

\noindent {\small Fig.4. Display of divisible by 3 natural
numbers` cobweb poset.}
\end{center}

\vspace{2mm}

\begin{center}

\includegraphics[width=75mm]{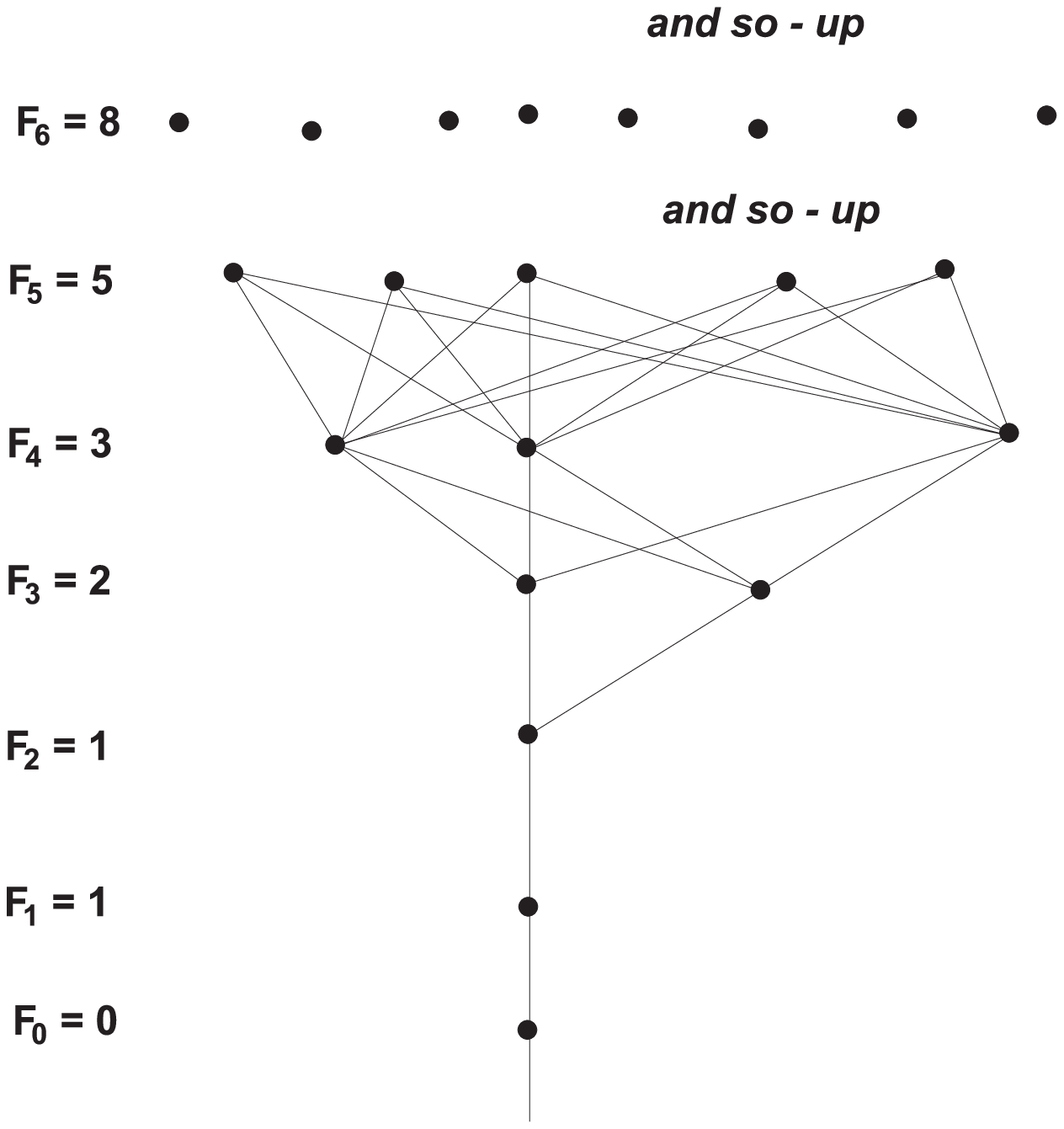}

\vspace{2mm}

\noindent {\small Fig.5. Display of Fibonacci numbers` cobweb
poset.}
\end{center}

Compare with  the bottom $6$ levels of a Young-\textit{Fibonacci
lattice}, introduced by Richard Stanley- in Curtis Greene`s
gallery of posets:
www.haverford.edu/math/cgreene/posets/posetgallery.html.

\vspace{2mm} \noindent As seen above - for example  the $Fig.5$.
displays the rule of the construction of the  Fibonacci "cobweb"
poset. It is being visualized clearly while defining this cobweb
poset  $P$ with help of its incidence matrix . The incidence
$\zeta$ function  matrix representing uniquely just this cobweb
poset $P$ has the staircase structure correspondent with
$"cobwebed"$ Fibonacci Tree i.e. a Hasse diagram of the particular
partial order relation under consideration. This is seen below on
the Fig.$6$ [8]:

\vspace{2mm}

$$ \left[\begin{array}{ccccccccccccccccc}
1 & 1 & 1 & 1 & 1 & 1 & 1 & 1 & 1 & 1 & 1 & 1 & 1 & 1 & 1 & 1 & \cdots\\
0 & 1 & 1 & 1 & 1 & 1 & 1 & 1 & 1 & 1 & 1 & 1 & 1 & 1 & 1 & 1 & \cdots\\
0 & 0 & 1 & 1 & 1 & 1 & 1 & 1 & 1 & 1 & 1 & 1 & 1 & 1 & 1 & 1 & \cdots\\
0 & 0 & 0 & 1 & 0 & 1 & 1 & 1 & 1 & 1 & 1 & 1 & 1 & 1 & 1 & 1 & \cdots\\
0 & 0 & 0 & 0 & 1 & 1 & 1 & 1 & 1 & 1 & 1 & 1 & 1 & 1 & 1 & 1 & \cdots\\
0 & 0 & 0 & 0 & 0 & 1 & 0 & 0 & 1 & 1 & 1 & 1 & 1 & 1 & 1 & 1 & \cdots\\
0 & 0 & 0 & 0 & 0 & 0 & 1 & 0 & 1 & 1 & 1 & 1 & 1 & 1 & 1 & 1 & \cdots\\
0 & 0 & 0 & 0 & 0 & 0 & 0 & 1 & 1 & 1 & 1 & 1 & 1 & 1 & 1 & 1 & \cdots\\
0 & 0 & 0 & 0 & 0 & 0 & 0 & 0 & 1 & 0 & 0 & 0 & 0 & 1 & 1 & 1 & \cdots\\
0 & 0 & 0 & 0 & 0 & 0 & 0 & 0 & 0 & 1 & 0 & 0 & 0 & 1 & 1 & 1 & \cdots\\
0 & 0 & 0 & 0 & 0 & 0 & 0 & 0 & 0 & 0 & 1 & 0 & 0 & 0 & 1 & 1 & \cdots\\
0 & 0 & 0 & 0 & 0 & 0 & 0 & 0 & 0 & 0 & 0 & 1 & 0 & 1 & 1 & 1 & \cdots\\
0 & 0 & 0 & 0 & 0 & 0 & 0 & 0 & 0 & 0 & 0 & 0 & 1 & 1 & 1 & 1 & \cdots\\
0 & 0 & 0 & 0 & 0 & 0 & 0 & 0 & 0 & 0 & 0 & 0 & 0 & 1 & 0 & 0 & \cdots\\
0 & 0 & 0 & 0 & 0 & 0 & 0 & 0 & 0 & 0 & 0 & 0 & 0 & 0 & 1 & 0 & \cdots\\
0 & 0 & 0 & 0 & 0 & 0 & 0 & 0 & 0 & 0 & 0 & 0 & 0 & 0 & 0 & 1 & \cdots\\
. & . & . & . & . & . & . & . & . & . & . & . & . & . & . & . & . \cdots\\
 \end{array}\right]$$

 \vspace{2mm}

\noindent \textbf{Figure 6.  The  staircase structure  of
incidence matrix $\zeta$ for the Fibonacci cobweb poset case}

\vspace{2mm}

\noindent \textbf{Note} The knowledge of $\zeta$  matrix explicit
form enables one  to construct (count) via standard algorithms
[10] the M{\"{o}}bius matrix $\mu =\zeta^{-1} $ and other typical
elements of incidence algebra perfectly suitable for calculating
number of chains, of maximal chains etc. in finite sub-posets of
$P$. All elements of the corresponding incidence algebra are then
given by a matrix of the Fig.5 with $1$`s replaced by any reals (
or ring elements in more general cases).

\vspace{2mm} \noindent Right from the definition of $P$ via its
Hasse diagram  here now follow quite obvious and important
observations. They lead us to a combinatorial interpretation of
cobweb poset`s characteristic binomial-like coefficients (for
example - fibonomial ones [4,5]).  Here they are with the first
obvious observation  at the start.

\vspace{3mm}
\noindent {\bf Observation 1}

{\it The number of maximal chains starting from The Root  (level
$0_F$) to reach any point at the $n-th$ level  with $n_F$ vertices
is equal to $n_{F}!$}.

\vspace{2mm}

\noindent {\bf Observation 2} $(k>0)$

{\it The number of maximal chains \textbf{rooted in any } vertex
at the $k-th$ level reaching the $n-th$ level with $n_F$ vertices
is equal to $n_{F}^{\underline{m}}$, \quad where $m+k=n.$ } \\

\vspace{2mm}\noindent Indeed. Denote the number of ways to get
along
maximal chains from  \textbf{any point} in $\Phi_k$ to $\Rightarrow  \Phi_n , n>k$ with the symbol\\
  $$[\Phi_k \rightarrow \Phi_n]$$
  then obviously we have :\\
           $$[\Phi_0 \rightarrow \Phi_n]= n_F!$$ and
$$[\Phi_0 \rightarrow \Phi_k]\times [\Phi_k\rightarrow \Phi_n]=
[\Phi_0 \rightarrow \Phi_n].$$

\vspace{2mm}

\noindent In order to formulate the combinatorial interpretation
of $F-sequence-nomial$ coefficients ({\it F-nomial} - in short)
[6,4,5,8] let us consider all finite \textit{"max-disjoint"}
sub-posets rooted at the $k-th$ level at any fixed vertex $\langle
r,k \rangle, 1 \leq r \leq k_F $  and ending  at corresponding
number of vertices at the $n-th$ level ($n=k+m$) where the
\textit{"max-disjoint"} sub-posets are defined below.

\vspace{2mm}

\begin{defn}
Two isomorphic copies of $P_m$ are said to be max-disjoint if
being considered as sets of maximal chains they are disjoint i.e
they have no maximal chain in common. All of $P_m$`s constitute
from now on a family of prime [1] $prefabiants$.
\end{defn}

\vspace{2mm}

\begin{defn}
We denote the number of all max-disjoint isomorphic copies of
$P_m$)  rooted at any vertex $\langle j,k \rangle , 1\leq j \leq
k_F $ of $k-th$ level  with the symbol
$$ \left( \begin{array}{c} n\\k\end{array}\right)_{F}.$$
\end{defn}
We use the accustomed to practical convention:  $\left(
\begin{array}{c} 0\\0\end{array}\right)_{F}=1.$

\noindent Naturally the above definition make sense not for
arbitrary $F$ sequences as  $F-nomial$ coefficients should be
nonnegative integers.

\vspace{2mm}

\begin{defn}
A sequence  $F = \{n_F\}_{n\geq 0}$ is called cobweb-admissible
iff
$$ \left( \begin{array}{c} n\\k\end{array}\right)_{F}\in N\cup\{0\}
\quad for \quad k,n\in N\cup\{0\}\equiv Z_\geq.$$
\end{defn}
Recall now that  the number of ways to reach an upper level from a
lower one along any of  maximal chains  i.e.  the number of all
maximal chains from the level
$\Phi_k $ to $ \Rightarrow  \Phi_n ,\quad n>k$ is equal to\\

  $$ [\Phi_k \rightarrow \Phi_n]= n_{F}^{\underline{m}}.$$

\noindent Naturally then we have

\begin{equation}
\left( \begin{array}{c} n\\k\end{array}\right)_{F} \times [\Phi_0
\rightarrow \Phi_m] = [\Phi_k \rightarrow \Phi_n]=
n_{F}^{\underline{m}}
\end{equation}
where  $[\Phi_0 \rightarrow \Phi_m]= m_F!$ counts the number of
maximal chains in any copy of the $P_m$. With this in mind we see
that the following holds.

\vspace{3mm}

\noindent {\bf Observation 3} $(\textbf{n,k}\geq\textbf{0})$

{\it Let $n = k+m$. The number of \textit{max-disjoint} sub-posets
isomorphic to $P_{m}$ (max-disjoint isomorphic copies of prime
$prefabiants$) , rooted at the $k-th$ level and ending at the n-th
level is equal to}

$$\frac{n_{F}^{\underline{m}}}{m_{F}!} =
\left( \begin{array}{c} n\\m\end{array}\right)_{F}$$
$$ = \left(\begin{array}{c} n\\k\end{array} \right)_{F}=
\frac{n_{F}^{\underline{k}}}{k_{F}!}. $$

\vspace{2mm} \noindent \textbf{Note} The Observation 3 provides us
with the \textbf{new combinatorial interpretation} of  the  class
\textbf{of all classical $F-nomial$ coefficients} including
distinguished binomial or distinguished Gauss $q$- binomial ones
or Konvalina generalized binomial coefficients of the first and of
the second kind [11]- which include Stirling numbers too. The vast
family of Ward-like [6] admissible by
$\psi=\langle\frac{1}{n_F!}\rangle_{n\geq 0}$-extensions
$F$-sequences [7,12] includes also those desired here which shall
be called \textit{ "GCD-morphic"} sequences. This means that
$GCD[F_n,F_m] = F_{GCD[n,m]}$ where $GCD$ stays for Greatest
Common Divisor operator. The Fibonacci sequence is  a much
nontrivial and guiding famous  example of GCD-morphic sequence.
Naturally  \textbf{incidence coefficients of any reduced incidence
algebra of full binomial type } [13] are GCD-morphic sequences
therefore they are now independently given a new cobweb
combinatorial interpretation via Observation 3. More on that - see
the next section where prefab combinatorial description is being
served. Before that - on the way - let us formulate  the following
problem (open?).

\vspace{2mm}

\noindent \textbf{Problem 1} \textit{Find effective
characterizations of the cobweb admissible sequence i.e.  find all
examples. }

\section{Cobweb posets as prefabs with $nonassociative$ and noncommutative composition}

Finite  cobweb sub-posets i.e. isomorphic copies of $P_m,  m\geq
0$ constitute connected acyclic digraphs as well as the Hasse
diagram of the infinite cobweb poset $P$ is. \textbf{D}irected
\textbf{a}cyclic \textbf{g}raphs are denoted as DAG`s [9]. Hence
one might call connected DAG`s - directed trees. As for the recent
development on acyclic digraphs we refer to [14] and references
therein. As in [14] one considers here digraphs on labeled
vertices and a "digraph" means a simple graph with at most one
edge directed from vertex to vertex. Loops and cycles of length
two are permitted in general, but parallel edges are forbidden.
"\textbf{Acyclic}" means that there are no cycles of any length.
Apart from Theorem 1 there note in [14] also Bibliographic remarks
on acyclic digraphs refereing to Robinson and  Stanley and then to
Bender et al. and Gessel. Because of an easy access to Plotnikov`s
paper [9] we shall take other definitions from there - if needed -
for granted. These are temporarily used just for the guiding
observation relating cobweb prefabs` digraphs to [9]. Namely - in
terminology of [9] - we make rather obvious observation.

\vspace{2mm}

\noindent {\bf Observation 4} {\it The Hasse (here upward
oriented) diagram of any  prime cobweb poset or $P$ is an oDAG.}

\vspace{2mm}

\noindent For the sake of explanation we quote after [9]: A poset
$P$ is of the {\it dimension} 2;  $dim\ P=2$ if there exist two
chains $L_1$ and $L_2$ such that $P = L_{1}\cap L_{2}$.

\noindent A digraph $G$ is called {\it the $\textbf{o}rderable$
digraph (oDAG)} if there exists a $dim\ 2$ poset such that its
Hasse diagram coincides with the digraph $G$.

\vspace{1mm}

\noindent We shall pass over now to the brief presentations of a
cobweb prefab combinatorial  structure [1] in which each object
($prefabiant$) is uniquely representable by construction as a
synthesis (composition) of powers of prime objects  where here
these are the cobweb sub-posets $P_n$ of $P$ that are to be
identified with prime $prefabiants$. Since now on we shall adhere
to the notation and terminology of [1]. We assume the acquaintance
of [1] which is justly considered as famous as important.

\vspace{1mm}

\noindent The definition of prefab combinatorial structure
$(S,\odot,f)$ here is assumed to be given by Definition 1 from [1]
except for associativity requirement $a_1)$  and commutativity
requirement $a_2)$, which are postponed until stated otherwise. In
general
 $a\odot b\neq b\odot a , a,b \in S$ - already for prime objects.
The definition of weighted (not necessarily associative,
commutative) prefab and enumerator $g(A), A\subseteq S$ are then
Definitions 2 an 3 from [1] correspondingly. We shall now
formulate Observation 5 (to to be checked by careful examination)-
observation  of distinguished importance for the combinatorial
interpretation of the \textbf{property $c_2)$ from the Definition
1} [1] of the prefab. The property $c_2)$  postulate from [1] is

\begin{equation}
|a\odot b|= \frac{f(a\odot b)}{f(a)f(b)},\quad  a,b\in S ,
\end{equation}
whenever   $a$ and $b$ have no common factor different from
identity prefabiant $i$ [1] where here $|A|$ denotes here  the
number of max-disjoint isomorphic copies of prime prefabiants in
the set $A=a\odot b$. The function $f$ satisfies the requirement
$c_1$ - of course.

\vspace{1mm}

\noindent {\bf Observation 5} {\it Let the enumerator or
generating function for prefab subsets be defined as indicated
above. Then the set of requirements $Prefab_{c(2)}= \{a_3), b_1),
b_2), c_1), c_2)\}$ is equivalent to set of requirements
$Prefab_{(Th.1)}= \{a_3), b_1), b_2), c_1), Theorem.1 \}$ , where
Theorem 1 means} Theorem 1 from [1].

\vspace{1mm}

\noindent Both sets of requirements define on $S$ the same prefab
structure (not necessarily commutative and associative ) where
requirements $b_1), b_2)$ are to be understood as rewritten in an
order and brackets being taken into account fashion.

\vspace{1mm}

\noindent Now comes the example of the class of weighted prefabs
$(S,\odot,f,\omega)$ with noncommutative, $nonassociative $
synthesis (composition) $\odot$. We shall call this binary
multivalued operation` [1] analogue case here a
"\textit{coopt-synthesis}" $\odot$. \vspace{2mm}

\noindent {\bf Cobweb prefab combinatorial structure.} The family
$S$ of combinatorial objects ($prefabiants$) consists of all
layers   $\langle\Phi_k \rightarrow \Phi_n \rangle,\quad k<n,\quad
k,n \in N\cup\{0\}\equiv Z_\geq$ and an empty prefabiant $i$.
Layer is considered here to be the set of  all max-disjoint
isomorphic copies (iso-copies) of $P_{n-k}$. The set $\wp$ of
prime objects consists of all sub-posets $\langle\Phi_0
\rightarrow \Phi_m \rangle$  i.e. all $P_m$`s $m \in
N\cup\{0\}\equiv Z_\geq$ constitute from now on a family of prime
$prefabiants$ [1]. The $Z_\geq$ grading preserving $\odot$
coopt-synthesis for prime prefabiants  $P_k \odot P_m$ =
$\langle\Phi_k \rightarrow \Phi_n\rangle, n= k+m$ means: consider
the leafs of  $P_k$  to be the roots of all max-disjoint
isomorphic copies (iso-copies) of $P_m$. Run through all the leafs
(now - roots). The $Z_\geq$ grading preserving $\odot$ synthesis
of (not necessarily prime) prefabiants - accordingly means the
same procedure with the requirement added (see: Example 5 in [1]).
If this algorithm applied to subsequent prime elements of the
second prefab gives rise to a layer of max-disjoint prefabs more
then one way - keep only one copy of it. As a result we have:

$$\langle\Phi_k \rightarrow \Phi_n \rangle \odot P_s = \langle\Phi_n \rightarrow
\Phi_{n+s}\rangle, \quad k\in Z_\geq , s \in N, n>k.$$

\noindent Accordingly the $Z_\geq \times Z_\geq $ grading of $S$
preserving $\odot$ synthesis ($ \odot$ \textit{coopt-synthesis})
is defined for arbitrary elements of $S$ as simply as follows:

$$\langle\Phi_k \rightarrow \Phi_n\rangle\odot \langle\Phi_t \rightarrow \Phi_{t+s}\rangle = \langle\Phi_n \rightarrow
\Phi_{n+s}\rangle;\quad t,k\in Z_\geq, n>k,s>0.$$

\vspace{1mm}

\noindent In order to satisfy the requirement $a_3)$ we postulate
for an empty prefabiant $i$ that
$$\langle\Phi_k \rightarrow \Phi_n \rangle \odot i = i \odot \langle\Phi_k \rightarrow
\Phi_n\rangle = \langle\Phi_k \rightarrow \Phi_n\rangle, \quad
k\in Z_\geq , n>k.$$

\noindent The appropriately adjusted requirements $b_1), b_2)$ are
satisfied by construction as

$$\langle\Phi_k \rightarrow \Phi_n \rangle = P_k \odot P_{n-k}, \quad k\in Z_\geq ,n>k.$$

\noindent The coopt-synthesis $\odot$ is nonassociative  by
construction as

$$(\langle\Phi_k \rightarrow \Phi_n\rangle\odot \langle\Phi_t \rightarrow \Phi_t+s\rangle)\odot\langle\Phi_p \rightarrow
\Phi_{p+q}\rangle = P_{n+s} \odot P_q ;\quad p,t,k\in Z_\geq,
n>k,s>0,q>0$$
while
$$\langle\Phi_k \rightarrow \Phi_n\rangle\odot (\langle\Phi_t \rightarrow \Phi_t+s\rangle\odot\langle\Phi_p \rightarrow
\Phi_{p+q}\rangle) = P_n \odot P_q ;\quad p,t,k\in Z_\geq,
n>k,s>0,q>0 . $$
As for the size functions let us aid with an
analogy (see: Example 5 in [1]).

\vspace{2mm}

\noindent \textbf{Analogy:}

$$Graphs ................................Cobweb-prefabs \quad\textbf{I}..............Cobweb-prefabs \quad\textbf{II}$$

$$vertices ......................................max-chains.............................  leafs\quad of\quad P_n`s$$

$$connected \quad G_n \quad on \quad [n]\quad\quad........P_n\quad cobweb ........................... P_n\quad cobweb $$

$$ size (G_n)= n  .................................size_1(P_n)=n! .......................size_2(P_n)= n$$

$$ f(G_n)=n!.......................................f(P_n)=n! ................................f(P_n)= n!\quad.$$

\vspace{1mm}

\noindent Recall now the Equation 3. We may draw now from all the
above the following \textbf{conclusion}.

\vspace{1mm}

\noindent \textbf{Conclusion I}

\noindent In the finite cobweb posets setting the $f$ function may
be chosen so as to be the $size_1$ of a prime prefabiant  with
$n_F$ leafs or so as to be factorial of the $size_2$ of a prime
prefab with $n_F$ leafs. This gives:
$$f(P_n)= n_F! , \quad f(P_m^k)= (km)_F!\quad$$
and the Equation 1 gets the required, expected combinatorial
interpretation for any cobweb prefab structure determined by the
choice of any sequence of natural numbers from the countless
family of cobweb admissible sequences. Thus we are equipped with
the cobweb prefab`s  uniform combinatorial interpretation of all
of them at once.

\vspace{1mm}

\noindent \textbf{Conclusion II}

\noindent Naturally  the Corollary 1  from [1] also holds in our
case. Choosing now the weight function to be of the form
$$\omega (a)= x^n,\quad n =size_2(a),\quad a\in S$$
we have the weighted  cobweb prefab and consequently (see:
Examples 5,10 in [1]) the formula for the cobweb weighted prefab
enumerator reads

\begin{equation}
g(S)= exp\{g(\wp)\},
\end{equation}
where
\begin{equation}
g(\wp) = exp_F\{x\}-1,
\end{equation}
while
\begin{equation}
exp_F\{x\}=\sum_{n\geq 0}\frac{x^{n}}{n_{F}!}
\end{equation}
$exp_F$  function [6,12,7,15] -  is the primary object of extended
finite operator calculus being recently developed in [7,12,15].
There it serves to define a central object of extended umbral
calculus i.e.  the  generalized translation operator
$E^{a}(\partial_{F})$ where the linear difference operator
$\partial_{F};\;\;\partial_{F}x^{n}=n_{F}x^{n-1};\;\;n\geq 0$ is
known under the name of the $F$-derivative  [6,7,14,15]. here
comes the example (11)  from [1] interpreted in the language of
$\psi$ - extensions [6] in their operator form [12,15,7].\\

\vspace{1mm}

\noindent \textbf{ Bender - Goldman - prefab example} Let the
"prefabian" $\hat q$-Bell numbers  $ B^{pref}_n( \gamma )$ be
defined as sums over $k$ of  $\hat S_q(n,k) $ Stirling numbers
equal  of the number of  unordered direct sums decompositions of
the $n$-dimensional vector space $V_{q,n}$ over $GF(q) \equiv F_q$
with $k$ summands. Then the Bender-Goldman exponential formula
(17) from [1] in $\psi$-extensions` notation [12,15,7] reads
$$ B^{pref}_{\gamma}(x)=
\sum_{n\geq
0}B^{pref}_n(\gamma)\frac{x^n}{n_{\gamma}!}=exp\{exp_{\gamma}(x) -
1\}. \quad\quad\quad(\gamma- e.g.f.)
$$
Here
$$ n_{\gamma}! = (q^n - 1)(q^n - q^1)...(q^n - q^{n-1})=
|GL_n(F_q)|,$$  $D_o(q)=1$ by convention  while $D_n(q) \equiv
B^{pref}_n(\gamma)=$ number of all unordered direct sums
decompositions of the vector space $V_{q,n}.$

\vspace{1mm}

\noindent \textbf{The natural hint}

\noindent The appealing analogy of the above schema and example
just presented give rise to questions on their eventual
correspondents  as $Stirling_F$ numbers of the second kind,
$exp_F$-ponential polynomials  and $F$-Dobinski like formulas.
Such extensions are more or less implicit in some papers . for
example - see Wagner`s (1.15) formula in [16] which formally
becomes of the $(\gamma- e.g.f.)$ formula form from above with now
almost arbitrary $\gamma =
\langle\frac{1}{n_{\gamma}!}\rangle_{n\geq 0}$  sequence (see also
[17] and references therein). These questions are to be considered
elsewhere. As for the related (determined by $F-nomial$`s)
extended umbral calculi in its operator form one may contact also
very recent review [18].

\vspace{2mm}

\section{Cobweb posets as prefabs with associative and commutative composition}

Here  another  single valued commutative and associative
composition case is presented in brief. The definition of the next
prefab combinatorial structure  with the single valued composition
$(S,\circ,f)$ is assumed to given  here by the Definition 1 from
[1] including  associativity requirement $a_1)$ and commutativity
requirement $a_2)$.

\noindent The family $S$ of combinatorial objects ($prefabiants$)
consists now  of all layers $\langle\Phi_k \rightarrow \Phi_n
\rangle, \quad k<n,\quad k,n \in N\cup\{0\}\equiv Z_\geq$ and an
empty prefabiant $i$ to be interpreted as the name or representant
of all "empty layers" $\langle\Phi_m \rightarrow \Phi_m \rangle.$
Layer is considered here as the set of all max-disjoint isomorphic
copies (iso-copies) of $P_{n-k}$. The set $\wp$ of prime objects
consists of all sub-posets \quad $\langle\Phi_0 \rightarrow \Phi_m
\rangle$  i.e. all $P_m$`s $m \in N\cup\{0\}\equiv Z_\geq$
constitute from now on a family of prime $prefabiants$ [1].

\vspace{1mm}

\noindent The $Z_\geq$ grading preserving $\circ$ coopt-synthesis
for prime $prefabiants$  $P_k \circ P_m$ = $\langle\Phi_0
\rightarrow \Phi_{n+m}\rangle , \quad n= k+m $\quad means:
consider the leafs of \quad $P_k$ \quad to be the
\textit{transitory}  roots of all max-disjoint isomorphic copies
(iso-copies) of $P_m$. Run through all the leafs (now -transitory
roots).

\noindent The $Z_\geq \times Z_\geq $ grading of $S$ grading
preserving $\circ$ synthesis of (not necessarily prime)
prefabiants - accordingly means the same procedure with the
requirement added (see: Example 5  in [1] ):  if this algorithm
applied to subsequent prime elements of the second prefab gives
rise to a layer of max-disjoint prefabs more then one way - keep
only one copy of it. As a result we have:

$$\langle\Phi_k \rightarrow \Phi_n \rangle \circ P_s = \langle\Phi_{k+0} \rightarrow
\Phi_{n+s}\rangle, \quad k\in Z_\geq , s \in N, n>k.$$

\noindent Accordingly the $Z_\geq \times Z_\geq $ grading of $S$
preserving \textit{coopt-synthesis} "$\circ$"  is defined as
follows:

$$\langle\Phi_k \rightarrow \Phi_n\rangle\circ\langle\Phi_p \rightarrow \Phi_q\rangle = \langle\Phi_{k+p} \rightarrow
\Phi_{n+q}\rangle;\quad p,k\in Z_\geq, n>k,q>p.$$

\vspace{2mm}

\noindent \textbf{Conclusion III}

\noindent In this  finite cobweb posets setting with associative
and commutative composition $\circ$ being single valued the $f$
function may be chosen as constant equal to $1$ function  (note
the other possibilities: $f(\langle\Phi_k \rightarrow
\Phi_n\rangle)= \alpha^{n-k},\quad \alpha\neq 0, n>k$).

\vspace{1mm}

\noindent This implies validity of the Corollary 2 and the
Corollary 3 from [1] also  in the cobweb "$\circ$ - case" -
providing us with direct efficient analogy to the cases of
unlabeled graphs, unordered partitions or factorizations of
integers (see Examples 1,2,3 in [1]).

\vspace{2mm}

\noindent  \textbf{Acknowledgements}

\noindent Discussions with Participants of Gian-Carlo Rota Polish
Seminar\\
$http://ii.uwb.edu.pl/akk/index.html$  - are highly appreciated.

\begin
{thebibliography}{99}
\parskip 0pt

\bibitem{1}
E. Bender, J. Goldman   {\it Enumerative uses of generating
functions} , Indiana Univ. Math.J. {\bf 20} 1971), 753-765.

\bibitem{2}
D. Foata and M. Sch"utzenberger, Th'eorie g'eometrique des
polynomes euleriens, (Lecture Notes in Math., No. 138).
Springer-Verlag, Berlin and New York, 1970.

\bibitem{3}
A. Nijenhuis and H. S. Wilf, Combinatorial Algorithms, 2nd ed.,
Academic Press, New York, 1978.

\bibitem{4}
A. K. Kwasniewski {\it Information on combinatorial interpretation
of Fibonomial coefficients }   Bull. Soc. Sci. Lett. Lodz Ser.
Rech. Deform. 53, Ser. Rech.Deform. {\bf 42} (2003), 39-41. ArXiv:
math.CO/0402291   v1 18 Feb 2004

\bibitem{5}
A. K. Kwa\'sniewski, {\it The logarithmic Fib-binomial formula}
Advanced Stud. Contemp. Math. {\bf 9} No 1 (2004), 19-26

\bibitem{6}
M. Ward: {\em A calculus of sequences}, Amer.J.Math. Vol.58,
(1936), 255-266.

\bibitem{7}
A. K. Kwa\'sniewski, {\it On simple characterizations of Sheffer
$\Psi$-polynomials and related propositions of the calculus of
sequences}, Bull.  Soc. Sci.  Lettres  \L \'od\'z {\bf 52},S\'er.
Rech. D\'eform. {\bf 36} (2002), 45-65. ArXiv: math.CO/0312397
$2003$

\bibitem{8}
A. K. Kwa\'sniewski, {\it More on combinatorial interpretation of
fibonomial coefficients} , Bull.  Soc. Sci.  Lettres  \L \'od\'z
{\bf 54},S\'er. Rech. D\'eform. {\bf 44} (2004), 23-38 --65.
ArXiv: math.CO/0402344 v1 22 Feb 2004

\bibitem{9}
A.D. Plotnikov {\it A formal approach to the oDAG/POSET problem}
(2004) html://www.cumulativeinquiry.com/Problems/solut2.pdf
(submitted to publication - March 2005)

\bibitem{10}
E.Krot: {\it The first ascent into the Fibonacci Cob-web Poset}
(submitted to publication - December 2004), ArXiv: math.CO/0411007

\bibitem{11}
J. Konvalina , {\it A Unified Interpretation of the Binomial
Coefficients, the Stirling Numbers and the Gaussian Coefficients}
The American Mathematical Monthly {\bf 107}(2000), 901-910.

\bibitem{12}
A. K. Kwa\'sniewski {\it Main  theorems of extended finite
operator calculus} Integral Transforms and Special Functions, {\bf
14} No 6 (2003), 499-516.

\bibitem{13}
E. Spiegel, Ch. J. O`Donnell  {\it Incidence algebras}  Marcel
Dekker, Inc. Basel $1997$.

\bibitem{14}
Brendan D. McKay, Frederique E. Oggier, Gordon F. Royle, N. J. A.
Sloane, Ian M. Wanless and Herbert S. Wilf, {\ Acyclic digraphs
and eigenvalues of (0,1)-matrices} Journal of Integer Sequences,
{\bf 7}, August 2004  Article 04.3.3  (arXiv: math.CO/0310423)

\bibitem{15}
A. K. Kwa\'sniewski: {\em On Extended Finite Operator Calculus of
Rota and Quantum Groups}, Integral Transforms and Special
Functions Vol 2, No 4, (2001), 333-340

\bibitem{16}
Carl G. Wagner {\it Generalized Stirling and Lah numbers} Discrete
Mathematics  {\bf 160} (1996), 199-218.

\bibitem{17}
A. K. Kwa\'sniewski: {\em Information on some recent applications
of umbral extensions to discrete mathematics} to appear in Review
Bulletin of Calcutta Mathematical Society Vol 13 ( 2005) ArXiv:
math.CO/0411145 7 Nov 2004

\bibitem{18}
A.K. Kwa\'sniewski, E. Borak: {\em Extended finite operator
calculus - an example of algebraization of analysis}  Central
European Journal of Mathematics {\bf 2} (5), 2005, 767-792. ArXiv
math.CO/0412233 14 Dec 2004

\end{thebibliography}



\end{document}